\documentclass{amsart}
\usepackage[all]{xy}
\usepackage{amssymb}

    \newtheorem{theorem}{Theorem}
\numberwithin{theorem}{section} \theoremstyle{plain}

\newtheorem{corollary}[theorem]{Corollary}

\newtheorem{proposition}[theorem]{Proposition}

\theoremstyle{definition}

\numberwithin{equation}{section}

\begin{document}
\title[Galois modules and fundamental groups]{A note on Galois modules and the algebraic fundamental group of
projective curves}
\author{Am\'{\i}lcar Pacheco}
\date{\today}
\address{Universidade Federal do Rio de Janeiro\\ Departamento de
Ma\-te\-m\'a\-ti\-ca Pura\\
Rua Guai\-aquil 83, Cachambi, 20785-050 Rio de Janeiro, RJ, Brasil}
\thanks{This work was partially supported by CNPq research grant 300896/91-3
and Pronex \#41.96.0830.00} \email{amilcar@impa.br}

\begin{abstract}
Let $X$ be a smooth projective connected curve of genus $g\ge 2$ defined over an algebraically closed field $k$ of
characteristic $p>0$. Let $G$ be a finite group, $P$ a Sylow $p$-subgroup of $G$ and $N_G(P)$ its normalizer in $G$. We
show that if there exists an \'etale Galois cover $Y\to X$ with group $N_G(P)$, then $G$ is the Galois group wan
\'etale Galois cover $\mathcal{Y}\to\mathcal{X}$, where the genus of $\mathcal{X}$ depends on the order of $G$, the
number of Sylow $p$-subgroups of $G$ and $g$. Suppose that $G$ is an extension of a group $H$ of order prime to $p$ by
a $p$-group $P$ and $X$ is defined over a finite field $\mathbb{F}_q$ large enough to contain the $|H|$-th roots of
unity. We show that integral idempotent relations in the group ring $\mathbb{C}[H]$ imply similar relations among the
corresponding generalized Hasse-Witt invariants.
\end{abstract}

\subjclass{14H30, 14G15}\maketitle

\section{Introduction}

Let $X$ be a smooth projective connected curve of genus $g\ge 2$ defined over an algebraically closed field $k$ of
characteristic $p>0$. Given a finite group $G$ when does there exist an \'etale Galois cover $Y\to X$ with group $G$?

Grothendieck gave a necessary and sufficient condition for this to take place when $|G|$ is prime to $p$, \cite{grot}.
Hasse-Witt theory does the same for a $p$-group $G$, \cite{ser1}. Later, the case of a group $G$ which is the extension
of a group $H$ of order prime to $p$ by a $p$-group $P$ necessary and sufficient conditions for $G$ to be such a Galois
group were obtained in \cite{pacste}. This result was extended for a group $G$ which is the extension of any group $H$
by a $p$-group $P$, \cite{bor}. The result in \cite{pacste} relied on the naive representation theoretic approach of
\cite{pac} as well as on techniques of embedding problems and cohomological dimension. The result in \cite{bor} relied
on modular representation theory.

One natural question is: how far can this type of approach give information on the construction \'etale Galois covers?
We start by discussing why it was needed in the two latter cases to take an extension of a group by a $p$-group. Next
we consider any finite group $G$ and take $P$ to be a Sylow $p$-subgroup of $G$ and let $N_G(P)$ be its normalizer in
$G$. Subject to the hypothesis that there exists an \'etale Galois cover $Y\to X$ with group $N_G(P)$, we show that $G$
is the Galois group of an \'etale Galois cover $\mathcal{Y}\to\mathcal{X}$ (using patching theory \cite{ste}), where
the genus of $\mathcal{X}$ depends on the order of $G$, the number of Sylow $p$-subgroups of $G$ and $g$.

Finally, the necessary and sufficient condition of \cite{pacste} can be expressed in terms of certain numbers, called
generalized Hasse-Witt invariants, which generalize the notion of the $p$-rank of a curve. Our next result is of a
different nature, however related to these numbers. At this point we assume that $H$ is a group of order prime to $p$
and the $X$ is defined over a finite field $\mathbb{F}_q$ large enough to contain the $|H|$-th roots of unity. We show
that integral idempotent relations in the group ring $\mathbb{C}[H]$ imply similar relations among the corresponding
generalized Hasse-Witt invariants. Analogous situations for the genera and Hasse-Witt invariants of curves were done
before in \cite{pac1} and \cite{kan}. As in \cite{pac1}, where relations among zeta-functions are first proved, the
result follows from relations between Artin $L$-functions.

\section{Galois groups of \'etale Galois covers}

Let $X$ be a smooth projective connected curve of genus $g\ge 2$ defined over an algebraically closed field $k$ of
characteristic $p>0$. Let $\pi_1(X)$ be its algebraic fundamental group. Let $\Gamma_g$ be the topological fundamental
group of a compact Riemann surface of genus $g$ and $\hat{\Gamma}_g$ its profinite completion. It is a result due to
Grothendieck \cite[Corollary 2.12]{grot} that there exists a surjective group homomorphism
$\varphi:\hat{\Gamma}_g\twoheadrightarrow\pi_1(X)$ whose kernel is contained in every open normal subgroup $N$ of
$\hat{\Gamma}_g$ such that $(\hat{\Gamma}_g:N)$ is prime to $p$. In particular, $\pi_1(X)$ is topologically finitely
generated. Let $\pi_A(X)$ be the set of isomorphism classes of finite groups which are quotients of $\pi_1(X)$, i.e.,
$G\in\pi_A(X)$ if and only if there exists an \'etale Galois cover $Y\to X$ such that $G=\text{Gal}(Y/X)$. Hence, the
determination of $\pi_1(X)$ is equivalent to that of $\pi_A(X)$.

Our first goal is to discuss $\pi_A(X)$. As a consequence of \cite[Corollary 2.12]{grot} Grothendieck obtained that a
finite group $G$ of order prime to $p$ lies in $\pi_A(X)$ if and only if $G$ is isomorphic to a finite quotient of
$\hat{\Gamma}_g$. The case of finite $p$-groups $P$ is first treated by Hasse-Witt theory, when $P$ is a $p$-elementary
abelian group \cite[\S11]{ser1}. In this case, $G\in\pi_A(X)$ if and only if it has at most $\gamma_X$ generators,
where $\gamma_X$ is the $\mathbb{F}_p$-dimension of the $p$-torsion subgroup $J_X[p]$ of the Jacobian variety $J_X$ of
$X$. This number is called the Hasse-Witt invariant or the $p$-rank of $X$. The general case is reduced to this one by
taking the quotient $P/\Phi(P)$ of $P$ by its Frattini subgroup $\Phi(P)=[P,P]P^p$ and using that $\pi_1(X)$ has
$p$-cohomological dimension $\text{cd}_p(\pi_1(X))$ at most $1$ \cite[Introduction]{pacste}.

The first mixed case to consider is that of a finite group $G$ which is an extension of a group $H$ of order prime to
$p$ by a $p$-group $P$ \cite[Theorem 1.3]{pacste}. We proved that the necessary and sufficient condition for
$G\in\pi_A(X)$ is that there exists an \'etale Galois cover $Y\to X$ such that $P/\Phi(P)$ injects in $J_Y[p]$ as an
$\mathbb{F}_p[H]$-module. Moreover, the \'etale Galois cover $Z\to X$ with group $G$ dominates $Y\to X$. In the next
section we will express this explicitly in terms of multiplicities of irreducible representations. It should be also
noted that in \cite[Corollary 3.5]{ste}, using formal patching, it was proved that a finite group $G$ having at most
$g$ generators lies in $\pi_A(\mathcal{X})$, for some smooth projective connected curve $\mathcal{X}$ of genus $g$. In
contrast, the first step of the proof of the former result is based in the case where $P$ is a $p$-elementary abelian
group and it was solved by naive representation theory in \cite[Propositions 2.4 and 2.5]{pac}. We will comment later
on how to go from this step to the general case.

Later Borne in \cite[Theorem 1.1]{bor} used modular representation theory to extend \cite[Theorem 1.3]{pacste} to the
case of a finite group $G$ which is the extension of any group $H$ by a $p$-group $P$. He also obtained a necessary and
sufficient condition for $G\in\pi_A(X)$ in terms of representation theory. At this point it is convenient to phrase
these two latter results in terms of embedding problems \cite[p. 366]{har1}. Let $\Pi$, $\mathcal{G}$, $\mathcal{H}$
and $\mathcal{P}$ be profinite groups. An embedding problem is a pair
$\mathcal{E}=(\mathcal{G}\twoheadrightarrow\mathcal{H},\Pi\twoheadrightarrow\mathcal{H})$ of surjective profinite group
homomorphisms. We say that $\mathcal{E}$ has a weak, respectively strong, solution if there exists a group
homomorphism, respectively a surjective group homomorphism, $\Pi\to\mathcal{G}$ such that the following diagram
commutes
$$\xymatrix{\mathcal{G}\ar@{->>}[r]&\mathcal{H}\\ &\Pi\ar@{->>}[u]\ar[ul]}.$$
The two latter situations can be described as obtaining a strong solution to the embedding problem
$$\xymatrix{0\ar[r]&P\ar[r]&G\ar[r]&H\ar[r]&0\\&&&\pi_1(X)\ar@{->>}[u]\ar[ul]}.$$
The kernel $\mathcal{P}$ of $\mathcal{G}\twoheadrightarrow\mathcal{H}$ is called the kernel of the embedding problem
$\mathcal{E}$. In the case where $\mathcal{P}$ is a profinite $p$-group, we say that $\mathcal{E}$ is a $p$-embedding
problem \cite[Section 2]{har2}. Thus, \cite[Theorem 1.3]{pacste} and \cite[Theorem 1.1]{bor} give a necessary and
sufficient condition for finite $p$-embedding problems for the fundamental group of projective curves to be strongly
solved.

One natural question is: how far can the Galois module theory approach go in order to obtain information on $\pi_A(X)$?
Let us first point out where the $p$-group kernel was used in both situations. In \cite[\S 6]{pacste}, we started from
the strong solution of the embedding problem $(\pi_1(X)\twoheadrightarrow H, P/\Phi(P)\rtimes H\twoheadrightarrow H)$.
This solution was obtained via naive representation theory, \cite[Propositions 2.4 and 2.5]{pac}. Then we used the
$p$-group kernel and $\text{cd}_p(\pi_1(X))\le 1$ to obtain via Galois cohomology \cite[Proposition 16, I-23]{ser2} a
weak solution to the embedding problem $(\pi_1(X)\twoheadrightarrow P/\Phi(P)\rtimes H,G\twoheadrightarrow
P/\Phi(P)\rtimes H)$. The strong solution is obtained via a Frattini type of argument, \cite[Remark 4.4]{pacste}. In
\cite[\S 2]{bor} the $p$-group kernel was used to obtain (via Clifford's theorem) a 1-1 correspondence between the set
of isomorphism classes of simple $k[H]$-modules and the set of isomorphism classes of simple $k[G]$-modules
\cite[Remark 2 to Proposition 2.4]{bor}. This was crucial in the proof of the reduction lemma \cite[Lemma 3.1]{bor},
which was used to prove \cite[Theorem 1.1]{bor}.

Abhyankar's conjecture is closely related to the problem of when a finite group $G$ lies in $\pi_A(X)$, \cite{abh}. We
now describe this conjecture. Let $x_1,\cdots,x_r\in X$ with $r\ge 1$ and $U=X-\{x_1,\cdots,x_r\}$. Let $\pi_1(U)$ be
the algebraic fundamental group of $U$ and $\pi_A(U)$ the set of isomorphism classes of finite groups which are
quotients of $\pi_1(U)$. The conjecture states that a finite group $G\in\pi_A(U)$ if and only if $G/p(G)\in\pi_A(U)$,
where $p(G)$ is the quasi-$p$ subgroup of $G$, i.e., the group generated by the Sylow $p$-subgroups of $G$. One side of
this conjecture was already known to Grothendieck. In fact, another consequence of \cite[Corollary 2.12]{grot} is that
a finite group $G$ of order prime to $p$ lies in $\pi_A(U)$ if and only if $G$ is a finite quotient of the profinite
completion $\hat{\Gamma}_{g,r}$ of the topological fundamental group $\Gamma_{g,r}$ of a compact Riemann surface of
genus $g$ minus $r$ points. If $G\in\pi_A(U)$, then $G/p(G)\in\pi_A(U)$ and Grothendieck's result give a necessary and
sufficient condition for this to take place. The converse was the issue in the conjecture. When $U$ is the affine line,
the conjecture was proved by Raynaud \cite{ray}, using rigid analytic techniques. Later, Harbater \cite{har3} extended
this result to any affine curse using formal geometry.

The analogue of Abhyankar's conjecture fails for projective curves. Indeed, suppose $X$ is an ordinary curve, i.e.,
$J_X[p]\cong(\mathbb{Z}/p\mathbb{Z})^g$. Let $n>g$ be an integer and $G=(\mathbb{Z}/p\mathbb{Z})^n$. Then
$G/p(G)=\{1\}\in\pi_A(X)$, however $G\notin\pi_A(X)$, since it violates the condition that the number of generators of
$G$ has to be at most $\gamma_X=g$. The representation theoretic conditions of \cite[Theorem 1.3]{pacste} and
\cite[Theorem 1.1]{bor} are the obstruction for the condition $G/p(G)\in\pi_A(X)$ to be sufficient for $G\in\pi_A(X)$.

The first idea to proceed is to first obtain a strong solution to the embedding problem $(\pi_1(X)\twoheadrightarrow
H,(p(G)/[p(G)])\rtimes H\twoheadrightarrow H)$ in a similar way to \cite[Propositions 2.4 and 2.5]{pac}, where $[p(G)]$
denotes the commutator of $p(G)$. Then we are lead to try to solve the embedding problem
$(\pi_1(X)\twoheadrightarrow(p(G)/[p(G)])\rtimes H,G\twoheadrightarrow(p(G)/[p(G)])\rtimes H)$. However, the kernel is
not any longer a $p$-group and the $\ell$-cohomological dimension of $\pi_1(X)$ is not necessarily $1$, for a prime
number $\ell\ne p$ (cf. \cite[Corollaire 4.3]{art} and \cite[\S 5]{pacste}). In the modular representation theory
approach, although \cite[Proposition 2.4]{bor} can be proved without the hypothesis of the embedding problem having a
$p$-kernel, \cite[Lemma 3.1]{bor} cannot. So this points to a bound on the Galois module theoretic method.

The strategy of strongly solving first a finite $p$-embedding problem for the affine line $\mathbb{A}^1$ in order to
show that $G/p(G)\in\pi_A(\mathbb{A}^1)$ implies $G\in\pi_A(\mathbb{A}^1)$ appears in \cite{ray}. There rigid patching
is used together with other techniques. Strongly solving finite $p$-embedding problems also appears as part of the
strategy of \cite{har3}, where also a prescribed local behavior is needed \cite{har2} and formal patching replaces
rigid patching. For an account of the sketch of both proofs see \cite[p. 85, 86]{har4}.

If we try to mimic the use of patching theory for projective curves, not only do we not realize $G$ as a Galois group
of an \'etale Galois cover of $X$, but also the base curve has its genus increased. This is due to \cite[Theorem
3.1]{ste} which we describe next.

Given an integer $g\ge 1$, let $\pi_A(g)$ be the set of isomorphism classes of finite groups $G$ such that
$G\in\pi_A(X)$, for some smooth projective connected curve $X$ of genus $g$ defined over $k$. Our main ingredient is
the following result.

\begin{proposition}\label{propkate}\cite[Theorem 3.1]{ste} Let $G$ be a finite group generated by two subgroups
$H_1$ and $H_2$. Suppose there exist smooth projective connected curves $X_1$ and $X_2$ of genus $g_1$ and $g_2$,
respectively, defined over $k$ such that $H_1\in\pi_A(U_1)$ and $H_2\in\pi_A(U_2)$, where $U_1=X_1-\{\xi_1\}$ and
$U_2=X_2-\{\xi_2\}$. Suppose furthermore that the ramification over $\xi_1$ and $\xi_2$ is tame and that the generator
of the inertia group of $\xi_1$ is the inverse in $G$ of that of the inertia group at $\xi_2$. Then
$G\in\pi_A(g_1+g_2)$.
\end{proposition}

Let $P$ be a Sylow $p$-subgroup of $G$ and $N_G(P)$ its normalizer in $G$. Let $|G|=mp^e$ with $p\nmid m$ and $n_p$ the
number of Sylow $p$-subgroups of $G$. By \cite[Chapter I, Theorem 3.7]{gor} $G=N_G(P)p(G)$. In particular, $G$ is
generated by $N_G(P)$ and $p(G)$.

\begin{theorem}Suppose there exists a smooth projective connected cur\-ve $X$ of genus $g\ge 2$ such that
$N_G(P)\in\pi_A(X)$ (hence satisfying the condition of \cite[Theorem 1.3]{pacste}). Then $G\in\pi_A(n_p+m(g-1)+g)$.
\end{theorem}

\begin{proof}Let $Z(P)\to X$ be an \'etale Galois cover with group $N_G(P)$. This cover dominates an \'etale
Galois cover $Y(P)\to X$ with group $N_G(P)/P$ and $P\in\pi_A(Y(P))$. By the Riemann-Hurwitz formula $Y(P)$ has genus
$(m/n_p)(g-1)+1$. Note that for any other Sylow $p$-subgroup $P'$ of $G$, $N_G(P')\cong N_G(P)$, hence
$N_G(P')\in\pi_A(X)$. Consequently, there exists an \'etale Galois cover $Z(P')\to X$ with group $N_G(P')$ which
dominates an \'etale Galois cover $Y(P')\to X$ with group $N_G(P')/P'$ and $P'\in\pi_A(Y(P'))$. By Proposition
\ref{propkate}, since $p(G)$ is generated by the Sylow $p$-subgroups of $G$, $p(G)\in\pi_A(m(g-1)+n_p)$. Let
$\mathcal{Y}'\to\mathcal{X}'$ be an \'etale Galois cover with group $p(G)$. Since $N_G(P)\in\pi_A(X)$, i.e., $N_G(P)$
is the Galois group of an \'etale Galois cover $Y\to X$, once again applying Proposition \ref{propkate}, we conclude
$G\in\pi_A(n_p+m(g-1)+g)$.
\end{proof}

\section{Idempotent relations and generalized $p$-ranks of curves}

Let $G$ be a finite group which is the extension of a group $H$ of order prime to $p$ by a $p$-group $P$. Let
$\Phi(P)=[P,P]P^p$ be the Frattini subgroup of $P$. Let $X$ be a smooth projective connected curve of genus $g\ge 2$
over $k$. The necessary and sufficient condition for $G$ to lie in $\pi_A(X)$ can be phrased in an explicit way using
representation theory.

For each irreducible character $\chi$ of $H$, let $V_{\chi}$ be the corresponding irreducible $k[H]$-module and
$\mathbb{F}_{p^{l_{\chi}}}=\mathbb{F}_p(\chi)$ the field generated by $\mathbb{F}_p$ and the values of $\chi$. Let
$\rho_{\chi}$ be the irreducible representation corresponding to $V_{\chi}$. Then for each positive power $p^m$ of $p$
we have another irreducible representation $\rho_{\chi^{p^m}}$ of $G$ such that
$\rho_{\chi^{p^m}}(\tau)=\rho_{\chi}(\tau)^{p^m}$, for every $\tau\in H$. Denote by $V_{\chi^{p^m}}$ the irreducible
$k[G]$-module corresponding to $\rho_{\chi^{p^m}}$. The irreducible $\mathbb{F}_p[H]$-modules are given by
$V_{[\chi]}=V_{\chi}\oplus V_{\chi^p}\oplus\ldots\oplus V_{\chi^{p^{l_{\chi}-1}}}$. Let
$P/\Phi(P)=\bigoplus_{[\chi]}V_{[\chi]}^{m_{[\chi]}}$ be the decomposition of $P/\Phi(P)$ into irreducible
$\mathbb{F}_p[H]$-modules.

Let $Y\to X$ a Galois cover with group $H$. The group $H$ acts on $J_Y[p]$. Let
$J_Y[p]\otimes_{\mathbb{F}_p}k=\bigoplus_{\chi}V_{\chi}^{r_{\chi}}$ be the decomposition of
$J_Y[p]\otimes_{\mathbb{F}_p}k$ into irreducible $k[H]$-modules and denote $J_Y[p]_{\chi}=V_{\chi}^{r_{\chi}}$. This is
a $k$-vector space whose dimension $\gamma_{Y,\chi}$ is called the generalized Hasse-Witt invariant of $Y$ with respect
to $\chi$ \cite[\S 2]{ruc}. It follows from \cite[Proposition 2.8 and Lemma 2.11]{pac} that
$\gamma_{Y,\chi^{p^m}}=\gamma_{Y,\chi}$ for every positive power $p^m$ of $p$. Thus, we denote it by
$\gamma_{Y,[\chi]}$. The condition, $P/\Phi(P)$ injects as an $\mathbb{F}_p[H]$-module in $J_Y[p]$ is equivalent to
$m_{[\chi]}n_{[\chi]}\le\gamma_{Y,[\chi]}$ for any irreducible $\mathbb{F}_p[H]$-module $V_{[\chi]}$, where
$n_{[\chi]}$ denotes the dimension of any of the $k[H]$-irreducible modules $V_{\chi^{p^m}}$. Hence, our representation
theoretic condition is expressed in terms of the generalized Hasse-Witt invariants.

These invariants are also used to count the number of \'etale Galois covers of $X$ whose Galois group $G$ is the
extension of a group $H$ of order prime to $p$ by an elementary $p$-abelian group $P$ such that the action of $H$ on
$P$ is faithful and irreducible \cite[Theorem 1.8]{pac} (for previous results see also \cite{nak} and \cite{kat}). When
instead of an algebraically closed field $k$ we consider a finite $\mathbb{F}_q$ large enough to contain the $|H|$-th
roots of unity, the generalized Hasse-Witt invariants are obtained as the degrees of the $L$-functions
$L(t,\chi,Y/X)^{n_{\chi}}$ which we now describe \cite[Corollary 4.1]{ruc}.

We suppose from now on that $Y$ and $X$ are defined over a finite field $\mathbb{F}_q$ of characteristic $p$ large
enough to contain the $|H|$-th roots of unity. Let $F:Y\to Y$ be the geometric Frobenius morphism of $Y$ with respect
to $\mathbb{F}_q$. For each positive integer $\nu$ and for each $\tau\in H$, let $\Lambda(\tau\circ F^{\nu})$ be the
number of points of $Y$ fixed by $\tau\circ F^{\nu}$. The Artin $L$-function $L(t,\chi,Y/X)$ of an irreducible
character $\chi$ of $H$ is defined as
$$L(t,\chi,Y/X)=\exp\left(\sum_{\nu\ge 1}A_{\nu}(\chi)\frac{t^{\nu}}{\nu}\right),\text{ where }A_{\nu}(\chi)=\frac
1{|H|}\sum_{\tau\in H}\chi(\tau^{-1})\Lambda(\tau\circ F^{\nu}).$$

Let $\ell\ne p$ be a prime number and let $T_{\ell}(J_Y)=\varprojlim J_Y[\ell^n]$ be the $\ell$-adic Tate module of
$J_Y$. It is a result due to Weil \cite[p. 186, Theorem 6]{lan} that there exists a representation
$\text{End}(J_Y)\to\text{End}(T_{\ell}(J_Y))$ given by $\tau\mapsto\alpha_{\tau}^*$. Let $\vartheta$ be the restriction
of this representation to $H$. Denote also by $F$ the geometric Frobenius of $J_Y$ over $\mathbb{F}_q$ and let
$\alpha_{F^{\nu}}^*\in\text{End}(T_{\ell}(J_Y))$ be the element corresponding to $\nu$-th power $F^{\nu}$ of $F$. Then
$\Lambda(\tau\circ F^{\nu})=q^{\nu}+1-\text{Tr}(\alpha_{\tau^{-1}}^*\circ\alpha_{F^{\nu}}^*\,|\,T_{\ell}(J_Y))$. Hence,
$$L(t,\chi,Y/X)=\frac{P(t,\chi,Y/X)}{(1-qt)^{u_{\chi}}(1-t)^{u_{\chi}}},$$
where \begin{equation}\label{pol1}P(t,\chi,Y/X)=\exp\left(\sum_{\nu\ge 1}-\frac 1{|H|}\left(\sum_{\tau\in
H}\chi(\tau^{-1})\text{Tr}(\alpha_{\tau^{-1}}^*\circ\alpha_{F^{\nu}}^*)\right)\frac{t^{\nu}}{\nu}\right)\end{equation}
is a polynomial in $t$ and $u_{\chi}=1$, if $\chi$ is trivial, $u_{\chi}=0$, otherwise. Denote by $\chi^{-1}$ the
character of the representation $\rho_{\chi}(\tau)^{-1}$, for every $\tau\in H$.

For each subgroup $M$ of $H$ and each irreducible character $\chi_{M}$ of $M$, let
$e_{\chi_M}=\frac{\chi_M(1)}{|M|}\sum_{\sigma\in M}\chi_M(\sigma^{-1})\sigma\in\mathbb{C}[H]$ be the corresponding
idempotent. An integral idempotent relation is an equality of the form $\sum_Ms_{\chi_M}e_{\chi_M}=0$, where $M$ runs
through the subgroups of $H$, $s_{\chi_M}\in\mathbb{Z}$ and for each $M$ we choose an irreducible character $\chi_M$ of
$M$.

\begin{theorem}\label{idempthm}
An integral idempotent relation $\sum_Ms_{\chi_M}e_{\chi_M}=0$ in $\mathbb{C}[H]$ implies a relation
$$\prod_MP(t,\chi_M^{-1},Y/Y_M)^{\chi_M(1)s_{\chi_M}}=1,$$
where $Y_M$ is the quotient curve $Y/M$ and $Y\to Y_M$ the corresponding Galois cover.
\end{theorem}

\begin{proof}Let $\mathbb{C}_{\ell}$ be the completion of the algebraic closure $\overline{\mathbb{Q}}_{\ell}$ of
$\mathbb{Q}_{\ell}$ and identify the $|H|$-th roots of unity in $\mathbb{C}$ with those in $\mathbb{C}_{\ell}$. We will
consider the integral idempotent relation $\sum_Ms_{\chi_M}e_{\chi_M}=0$ inside $\mathbb{C}_{\ell}[H]$ instead of
$\mathbb{C}[H]$. Extend the representation $\vartheta$ to a representation
$\overline{\vartheta}:\mathbb{C}_{\ell}[H]\to\text{GL}(T_{\ell}(J_Y)\otimes_{\mathbb{Z}_{\ell}}\mathbb{C}_{\ell})$.
Then
$$0=\overline{\vartheta}\left(\sum_Ms_{\chi_M}e_{\chi_M}\right)=\sum_M\left(\frac{s_{\chi_M}\chi_M(1)}
{|M|}\right)\sum_{\sigma\in M}\chi_M(\sigma^{-1})\alpha_{\sigma}^*.$$ Taking the trace of following representation in
$T_{\ell}(J_Y)\otimes_{\mathbb{Z}_{\ell}}\mathbb{C}_{\ell}$, we obtain
$$0=\text{Tr}\left(\overline{\vartheta}\left(\sum_Ms_{\chi_M}e_{\chi_M}\right)\circ\alpha_{F^{\nu}}^*
\right)=\sum_M \frac{s_{\chi_M}}{|M|}\sum_{\sigma\in
M}\chi_M(\sigma^{-1})\text{Tr}(\alpha_{\tau}^*\circ\alpha_{F^{\nu}}^*).$$ The result now follows from (\ref{pol1}).
\end{proof}

The Artin $L$-function $L(t,Y/X,\chi)$ can be expressed $p$-adically as follows. For each $P\in J_Y$, let
$\mathcal{O}_{J_Y,P}$ be the local ring of $J_Y$ at $P$. For each integer $n\ge 1$, let $W_n(\mathcal{O}_{J_Y,P})$ be
the Witt ring of length $n$ with coefficients in $\mathcal{O}_{J_Y,P}$. Let $\mathcal{W}_n$ be the sheaf over $J_Y$
defined by these rings. The cohomology groups $H^m(J_Y,\mathcal{W}_n)$ form a projective limit. Let
$H^m(J_Y,\mathcal{W})=\varprojlim H^m(J_Y,\mathcal{W}_n)$. Let $J_Y^{\wedge}$ be the dual abelian variety of $J_Y$ and
$T_p(J_Y^{\wedge})$ its $p$-adic Tate module. Let $W(\overline{\mathbb{F}}_q)$ the ring of Witt vectors of
$\overline{\mathbb{F}}_q$ and
$\overline{T}_p(J_Y^{\wedge})=T_p(J_Y^{\wedge})\otimes_{\mathbb{Z}_p}W(\overline{\mathbb{F}}_q)$. Let
$L(J_Y)=H^1(J_Y,\mathcal{W})\oplus\overline{T}_p(J_Y^{\wedge})$. Note that identifying the $|H|$-th roots of unity in
$\mathbb{C}$ with those in $\mathbb{F}_q$, we can assume that the characters of $H$ take value in $\mathbb{F}_q$, hence
in $W(\overline{\mathbb{F}}_q)$. Let $L(J_Y)=\bigoplus_{\chi}\mathcal{V}_{\chi}^{m_{\chi}}$ be the decomposition of
$L(J_Y)$ into irreducible $W(\overline{\mathbb{F}}_q)[H]$-modules and let
$L(J_Y)_{\chi}=\mathcal{V}_{\chi}^{m_{\chi}}$. It is proved in \cite[Corollary 3.1]{ruc} that
$P(t,\chi,Y/X)^{\chi(1)}=\det(1-\alpha_{p,F}^*t\,|\,L(J_Y)_{\chi^{-1}})$, where $\alpha_{p,F}^*\in L(J_Y)$ is the
endomorphism corresponding to the geometric Frobenius $F$ of $J_Y$. In particular, $P(t,\chi,Y/X)$ has $p$-integral
coefficients.

Let $K=\mathbb{F}_q(Y)$, let $\Omega_K$ be the space of all differentials of $K$ and $H^0(Y,\Omega_Y)$ the subspace of
regular differentials. Let $u\in K$ be a separating variable and $\omega=fdu\in\Omega_K$. The Cartier operator is
defined by $\mathcal{C}(\omega)=(-d^{p-1}f/du^{p-1})^{/1/p}du$. It is a $1/p$-linear operator, i.e.,
$\mathcal{C}(a^p\omega)=a\omega$, for $a\in K$, and it acts on $H^0(Y,\Omega_Y)$ \cite[\S 10, p.39]{ser1}.

Let $H^0(Y,\Omega_Y)=\bigoplus_{\chi}H^0(Y,\Omega_Y)_{\chi}$ be the decomposition of $H^0(Y,\Omega_Y)$ into a sum of
irreducible $\mathbb{F}_q[H]$-modules, where $H^0(Y,\Omega_Y)_{\chi}=V_{\chi}^{m_{\chi}}$ and $V_{\chi}$ is the
irreducible $\mathbb{F}_q[H]$-module corresponding to $\chi$. Let $q=p^n$, then $\mathcal{C}^n$ acts on
$H^0(Y,\Omega_Y)_{\chi}$. Moreover,
$\gamma_{Y,\chi}=\dim_{\mathbb{F}_q}\ker(1-\mathcal{C}^n\,|\,H^0(Y,\Omega_Y)_{\chi})$ \cite[Proposition 2.3]{ruc}. It
is also shown that $\det(1-\alpha_{p,F}^*t\,|\,L(J_Y)_{\chi^{-1}})\pmod
p=\det(1-\mathcal{C}^nt\,|\,H^0(Y,\Omega_Y)_{\chi})$ \cite[Theorem 4.1]{ruc}. Consequently,
$\deg(\det(1-\alpha_{p,F}^*t\,|\,L(J_Y)_{\chi^{-1}})\pmod p)=\gamma_{Y,\chi}$ and we obtain the following corollary of
Theorem \ref{idempthm}.

\begin{corollary}An integral idempotent relation $\sum_Ms_{\chi_M}e_{\chi_M}=0$ in $\mathbb{C}[H]$ implies a relation
$$\sum_Ms_{\chi_M}\gamma_{Y,\chi_M^{-1}}=0.$$
\end{corollary}

\end{document}